\documentclass[12pt,oneside]{amsart}
\usepackage{geometry}                % See geometry.pdf to learn the layout options. There are lots.
\usepackage{graphicx}
\usepackage{amssymb}
\usepackage{epstopdf}
\usepackage{amsmath,amsthm,amscd,amssymb}
\usepackage{latexsym}
\usepackage[colorlinks,citecolor=red,pagebackref,hypertexnames=false]{hyperref}
\numberwithin{equation}{section}
\theoremstyle{plain}
\newtheorem{theorem}{Theorem}[section]
\newtheorem{lemma}[theorem]{Lemma}

\theoremstyle{definition}

\theoremstyle{remark}
\newtheorem{remark}[theorem]{Remark}

\newtheorem{case[theorem]}{Case}

\date{October 30, 2011}      
\author{Alex Iosevich, Mihalis Mourgoglou and Krystal Taylor}
\address{Department of Mathematics, University of Rochester, Rochester, NY}
\email{iosevich@math.rochester.edu}
\address{D\'{e}partement de Math\'{e}matiques\\ UMR 8628 Universit\'{e} Paris-Sud 11-CNRS \\  B\^atiment 425\\ Facult\'{e} des Sciences d'Orsay\\ Universit\'{e} Paris-Sud 11\\ F-91405 Orsay Cedex}
\email{mihalis.mourgoglou@math.u-psud.fr}
\address{Department of Mathematics, University of Rochester, Rochester, NY}
\email{taylor@math.rochester.edu}
\thanks{This work was partially supported by the NSF Grant DMS10-45404.}
\title{\parbox{14cm}{\centering{On the Mattila-Sjolin theorem for distance sets}}}
\begin{document}
\maketitle
\begin{abstract} We extend a result, due to Mattila and Sjolin, which says that if the Hausdorff dimension of a compact set $E \subset {\Bbb R}^d$, $d \ge 2$, is greater than $\frac{d+1}{2}$, then the distance set $\Delta(E)=\{|x-y|: x,y \in E \}$ contains an interval. We prove this result for distance sets $\Delta_B(E)=\{ {||x-y||}_B: x,y \in E \}$, where ${|| \cdot ||}_B$ is the metric induced by the norm defined by a symmetric bounded convex body $B$ with a smooth boundary and everywhere non-vanishing Gaussian curvature. We also obtain some detailed estimates pertaining to the Radon-Nikodym derivative of the distance measure. 
\end{abstract}  
\maketitle
%\tableofcontents
\section{Introduction} 
The classical Falconer distance conjecture, originated in 1985 (\cite{Fal86}) says that if the Hausdorff dimension of a compact subset of ${\Bbb R}^d$, $d \ge 2$, is greater than $\frac{d}{2}$, then the Lebesgue measure of the set of distances, $\Delta(E)=\{|x-y|: x,y \in E \}$ is positive. Falconer (\cite{Fal86}) proved the first result in this direction by showing that ${\mathcal L}^1(\Delta(E))>0$ if the Hausdorff dimension of $E$ is greater than $\frac{d+1}{2}$. See also \cite{Falc86} and \cite{M95} for a thorough description of the problem and related ideas. The best currently known results are due to Wolff in two dimensions, and to Erdogan (\cite{Erd05}) in dimensions three and greater. They prove that ${\mathcal L}^1(\Delta(E))>0$ if the Hausdorff dimension of $E$ is greater than $\frac{d}{2}+\frac{1}{3}$. 

In another direction, an important addition to this theory is due to Mattila and Sjolin (\cite{MS99}) who proved that if the Hausdorff dimension of $E$ is greater than $\frac{d+1}{2}$, then $\Delta(E)$ not only has positive Lebesgue measure, but also contains an interval. This is accomplished by showing that the natural measure on the distance set has a continuous density. It was previously shown by Mattila (\cite{M85}) that if the ambient dimension is two or three, then the density of the distance measure is not in general bounded if the Hausdorff dimension of the underlying set $E$ is smaller than $\frac{d+1}{2}$. In higher dimensions, this question is still open for the Euclidean metric, but has been resolved if the Euclidean metric is replaced by a metric induced by a norm defined by a suitably chosen paraboloid. See \cite{IS10}. 

In this paper we give an alternative proof of the Mattila-Sjolin result and extend it to more general distance sets $\Delta_B(E)=\{ {||x-y||}_B: x,y \in E \}$, where ${|| \cdot ||}_B$ is the norm generated by a symmetric bounded convex body $B$ with a smooth boundary and everywhere non-vanishing Gaussian curvature. 

Our main result is the following.  
\begin{theorem} \label{main} Let $E$ be a compact subset of ${\Bbb R}^d$, $d \ge 2$, with Hausdorff dimension, denoted by $s$, greater than $\frac{d+1}{2}$. Let $\mu$ be a Frostman measure on $E$. Let $\sigma$ denote the Lebesgue measure on $\partial B$. Define the distance measure $\nu$ by the relation 
$$ \int h(t) d\nu(t)=\int \int h({||x-y||}_B) d\mu(x) d\mu(y),$$ where ${|| \cdot ||}_B$ is the norm generated by a symmetric bounded convex body $B$ with a smooth boundary and everywhere non-vanishing Gaussian curvature. 
\begin{itemize} 
\item{i)} Then the measure $\nu$ is absolutely continuous with respect to the Lebesgue measure. 

\vskip.125in 

\item{ii)} We have 
$$ \frac{\nu((t-\epsilon, t+\epsilon))}{2\epsilon}=M(t)+R^{\epsilon}(t), $$ where 
$$ M(t)=\int {|\widehat{\mu}(\xi)|}^2 \widehat{\sigma}(t \xi) t^{d-1} d\xi$$ is the density of $\nu$ and 
$$ \sup_{0<\epsilon<\epsilon_0} |R^{\epsilon}(t)| \lesssim \epsilon_0^{s-\frac{d+1}{2}}.$$ 

\vskip.125in 

\item{iii)} Moreover, $M \in C^{\lfloor s-\frac{d+1}{2} \rfloor}(I)$ for any interval $I$ not containing the origin, where $\lfloor u \rfloor$ denotes the smallest integer greater than or equal to $u$. In particular, $M$ is continuous away from the origin if $s>\frac{d+1}{2}$ and therefore $\Delta_B(E)$ contains an interval in view of i). 

\vskip.125in 

\item{iv)} Suppose that $s>k+\alpha$, where $k$ is a non-negative integer and $0<\alpha<1$. Then the $kth$ derivative of the density function of $\nu$ is H\"older continuous of order $\alpha$.
\end{itemize} 
\end{theorem} 

\begin{remark} Metric properties of ${|| \cdot ||}_B$ are not used in the proof of Theorem \ref{main}. Let $\Gamma$ be a star shaped body in the sense that for every $\omega \in S^{d-1}$ there exists $1<r_0(\omega)<2$ such that $\{r \omega: 0 \leq r \leq r_0(\omega) \} \subset \Gamma$ and $\{r \omega: r>r_0 \} \cap \Gamma=\emptyset$. Define ${||x||}_{\Gamma}=\inf \{t>0: x \in t\Gamma \}$ and let $\Delta_{\Gamma}(E)=\{ {||x-y||}_{\Gamma}: x,y \in E \}$. Let $\sigma_{\Gamma}$ denote the Lebesgue measure on the boundary of $\Gamma$. Then if $|\widehat{\sigma}_{\Gamma}(\xi)| \lesssim {|\xi|}^{-\frac{d-1}{2}}$,  the conclusion of Theorem \ref{main} holds with the same exponents. \end{remark} 

\subsection{Sharpness of results:} As we note above, Mattila's construction (\cite{M85}) shows that if the Hausdorff dimension of $E$ is smaller than $\frac{d+1}{2}$, $d=2,3$, then the density of distance measure is not in general bounded in the case of the Euclidean metric. Moreover, Mattila construction can be easily extended to all metrics generated by a bounded convex body $B$ with a smooth boundary and non-vanishing Gaussian curvature. 

In dimensions four and higher, all we know at the moment (see the main result in \cite{IS10}) is that there exists a bounded convex body $B$ with a smooth boundary and non-vanishing curvature, such that the density of the distance measure is not in general bounded if the Hausdorff dimension of the underlying set $E$ is less than $\frac{d+1}{2}$. We do no know what happens when the Hausdorff dimension of $E$ equals $\frac{d+1}{2}$ in any dimension and for any smooth metric. 

It would be very interesting if any of these results actually depended on the underlying convex body $B$ in a non-trivial way. This would mean that smoothness and non-vanishing Gaussian curvature of the level set do not tell the whole story. There is some evidence that this may be the case. See, for example, \cite{IR07}, where connections between problems of Falconer type and distribution of lattice points in thin annuli are explored. 

If the Hausdorff dimension of $E$ is less than $\frac{d}{2}$, then the density of the distance measure, for any metric induced by a bounded convex body $B$ with a smooth boundary and non-vanishing curvature is not in general bounded by a construction due to Falconer (\cite{Fal86}). 

\subsection{Acknowledgements} The second named author holds a two-year Sophie Germain International post-doctoral scholarship in Fondation de Math\'{e}matiques Jacques Hadamard (FMJH) and would like to thank the faculty and staff of the Universit\'{e} Paris-Sud 11, Orsay for their hospitality. 

\vskip.25in 

\section{Proof of Theorem \ref{main}} 

\vskip.125in 

\subsection{Proof of items i) and ii)} The proof of item i) of Theorem \ref{main} is due to Falconer (\cite{Fal86}) and Mattila (\cite{M85}). This brings us to item ii). Recall that every compact set $E$ in ${\Bbb R}^d$, of Hausdorff dimension $s>0$ possesses a Frostman measure (see e.g. \cite{M95}, p. 112), which is a probability measure $\mu$ with the property that for every ball of radius $r^{-1}$, denoted by $B_{r^{-1}}$, 
$$ \mu(B_{r^{-1}}) \lessapprox r^{-s},$$ where here, and throughout, $X \lessapprox Y$, with the controlling parameter $r$ means that for every $\epsilon>0$ there exists $C_{\epsilon}>0$ such that $X \leq C_{\epsilon}r^{\epsilon}Y$. Let 

$$\nu^{\epsilon}(t)=\frac{\nu((t-\epsilon, t+\epsilon))}{2\epsilon}=\frac{1}{2 \epsilon} \mu \times \mu \{(x,y): t-\epsilon \leq {||x-y||}_B \leq t+\epsilon \}.$$ 

We shall prove that $\lim_{\epsilon \to 0} \nu^{\epsilon}(t)$ exists and is a $C^{\lfloor s-\frac{d+1}{2} \rfloor}$ function. 

\vskip.125in 

Let $\rho$ be a smooth cut-off function, identically equal to $1$ in the unit ball and vanishing outside the ball of radius $2$. Let $\rho_{\epsilon}(x)=\epsilon^{-d} \rho(x/\epsilon)$. It is not hard to check that one can construct a $\rho$ such that the difference between $\nu^{\epsilon}(t)$ and 
\begin{equation} \label{cooler} \int \int \sigma_t*\rho_{\epsilon}(x-y) d\mu(x) d\mu(y),\end{equation} where $\sigma_t$ is the surface measure on the set $\{x: {||x||}_B=t \}$, is $o(\epsilon)$, so there is no harm in taking (\ref{cooler}) as the definition of $\nu^{\epsilon}(t)$. By the Fourier inversion formula, 
$$ \nu^{\epsilon}(t)=\int {|\widehat{\mu}(\xi)|}^2 \widehat{\sigma}_t(\xi) \widehat{\rho}(\epsilon \xi) d\xi$$
$$=\int {|\widehat{\mu}(\xi)|}^2 \widehat{\sigma}_t(\xi) d\xi-\int {|\widehat{\mu}(\xi)|}^2 \widehat{\sigma}_t(\xi) (1-\widehat{\rho}(\epsilon \xi)) d\xi$$
$$=M(t)+R^{\epsilon}(t).$$ 

We shall prove that $M(t)$ is a $C^{\lfloor s-\frac{d+1}{2} \rfloor}$ function and that $\lim_{\epsilon \to 0} R^{\epsilon}(t)=0$. We start with the latter. We shall need the following stationary phase estimate. See, for example, \cite{St93}, \cite{So93} or \cite{W04}. 
\begin{lemma} \label{stationary} Let $\sigma$ be the surface measure on a compact piece of a smooth convex surface in ${\Bbb R}^d$, $d \ge 2$, with everywhere non-vanishing Gaussian curvature. Then 
$$ |\widehat{\sigma}(\xi)| \lesssim {|\xi|}^{-\frac{d-1}{2}},$$ where here, and throughout, $X \lesssim Y$ means that there exists $C>0$ such that $X \leq CY$. 
\end{lemma} 

We shall also need the following well-known estimate. See, for example, \cite{Falc86} and \cite{M95}. 
\begin{lemma} \label{energy} Let $\mu$ be a Frostman measure on a compact set $E$ of Hausdorff dimension $s>0$. Then 
$$ \int_{2^j \leq |\xi| \leq 2^{j+1}} {|\widehat{\mu}(\xi)|}^2 d\xi \lessapprox 2^{j(d-s)},$$ and, consequently, 
$$ \int {|\widehat{\mu}(\xi)|}^2 {|\xi|}^{-\gamma} d\xi=c\int \int {|x-y|}^{-d+\gamma} d\mu(x)d\mu(y) \lesssim 1$$ if $\gamma>d-s$. 
\end{lemma} 

To prove the lemma, observe that 
$$  \int_{2^j \leq |\xi| \leq 2^{j+1}} {|\widehat{\mu}(\xi)|}^2 d\xi \lesssim \int {|\widehat{\mu}(\xi)|}^2 \psi(2^{-j}\xi) d\xi,$$ where $\psi$ is a suitable smooth function supported in $\{x \in {\Bbb R}^d: 1/2 \leq |x| \leq 4 \}$ and identically equal to $1$ in the unit annulus. By definition of the Fourier transform and the Fourier inversion theorem, this expression is equal to 
$$ 2^{dj} \int \int \widehat{\psi}(2^j(x-y)) d\mu(x) d\mu(y) \lessapprox 2^{j(d-s)}$$ since $\widehat{\psi}$ decays rapidly at infinity. 

By Lemma \ref{stationary} and Lemma \ref{energy}, we have 
$$ |R^{\epsilon}(t)| \lessapprox \int_{|\xi|>\frac{1}{\epsilon}} 
{|\widehat{\mu}(\xi)|}^2 {|\xi|}^{-\frac{d-1}{2}} d\xi$$ 
$$ \leq \int_{|\xi|>\frac{1}{\epsilon}} {|\widehat{\mu}(\xi)|}^2 {|\xi|}^{-\frac{d-1}{2}} d\xi=\sum_{j>\log_2(1/\epsilon)} \int_{2^j \leq |\xi| \leq 2^{j-1}} {|\widehat{\mu}(\xi)|}^2 {|\xi|}^{-\frac{d-1}{2}} d\xi$$ 
\begin{equation} \label{boring} \lessapprox \sum_{j>\log_2(1/\epsilon)} 2^{j(d-s)} 2^{-j \frac{d-1}{2}} \lessapprox \epsilon^{s-\frac{d+1}{2}}, \end{equation} and thus $\lim_{\epsilon \to 0} R^{\epsilon}(t)=0$. (To handle $|R^{\epsilon}(t)|$ over the integral when ${|\xi|<\frac{1}{\epsilon}}$, we notice that $(1-\widehat{\rho}(\epsilon \xi))$ is $0$ when $\xi=0$ and, by continuity, is small in a neighborhood about $0$. Then, we may dilate and re-define $\rho$ to get that $(1-\widehat{\rho}(\epsilon \xi))$ is small in a neighborhood of our choosing.)  This calculation establishes all the claims in part ii) of Theorem \ref{main}  

\vskip.125in 

\subsection{Proof of item iii)} Once again, by Lemma \ref{stationary}, we have 
$$ |M(t)| \lesssim \int {|\widehat{\mu}(\xi)|}^2 {|\xi|}^{-\frac{d-1}{2}} d\xi$$ and by the calculation identical to the one in the previous paragraph, we see that this quantity is $\lesssim 1$ if the Hausdorff dimension of $E$ is greater than $\frac{d+1}{2}$. Continuity follows by the Lebesgue dominated convergence theorem. The convergence of the integral allows us to differentiate inside the integral sign. We obtain 
$$ M'(t)=\int {|\widehat{\mu}(\xi)|}^2 \frac{d}{dt} \left\{t^{d-1} \widehat{\sigma}(t \xi) \right\} d\xi.$$ 
We have 
$$  \frac{d}{dt} \left\{t^{d-1} \widehat{\sigma}(t \xi)\right\}
=(d-1)t^{d-2} \widehat{\sigma}(t \xi)+t^{d-1} \nabla \widehat{\sigma}(t \xi) \cdot \xi.$$ 
Applying Lemma \ref{stationary} once more, the best we can say is that 
$$  \left| \frac{d}{dt} \left\{t^{d-1} \widehat{\sigma}(t \xi)\right\} \right| \lesssim {|\xi|}^{-\frac{d-1}{2}+1}.$$ 
Repeating the argument in \ref{boring}, we see that $M'(t)$ exists if the Hausdorff dimension of $E$ is greater than $\frac{d+1}{2}+1$. Proceeding in the same way one establishes that 
$$ \frac{d^m}{dt^m} \left\{ t^{d-1} \widehat{\sigma}(t \xi) \right\} \lesssim {|\xi|}^{-\frac{d-1}{2}
+m}$$ and the conclusion of Theorem \ref{main} follows. 
\subsection{Proof of item iv)} We shall deal with the case $k=0$, as the other cases follow from a similar argument. Let 
$$\lambda(t)=t^{d-1}\widehat{\sigma}(t \xi).$$ 

We must show that 
$$ |M(u)-M(v)| \leq C{|u-v|}^{\alpha}.$$ 

We have 
$$ M(u)-M(v)=\int {|\widehat{\mu}(\xi)|}^2 (\lambda(u)-\lambda(v)) d\xi$$
$$=\int {|\widehat{\mu}(\xi)|}^2 {(\lambda(u)-\lambda(v))}^{\alpha} {(\lambda(u)-\lambda(v))}^{1-\alpha} d\xi.$$

Now, 
$$ \lambda(u)-\lambda(v)=(u-v) \lambda'(c),$$ where $c \in (u,v)$, by the mean-value theorem. It follows that 
$$ {|\lambda(u)-\lambda(v)|}^{\alpha} \leq {|u-v|}^{\alpha} {|\lambda'(c)|}^{\alpha}.$$ 

On the other hand, 
$$ {|\lambda(u)-\lambda(v)|}^{1-\alpha} \leq {|\lambda(u)|}^{1-\alpha}+{|\lambda(v)|}^{1-\alpha}.$$ 

We have already shown above that 
$$ |\lambda(u)| \lesssim {|\xi|}^{-\frac{d-1}{2}} \ \text{and} \ |\lambda'(u)| \lesssim {|\xi|}^{-\frac{d-1}{2}+1}.$$ 

It follows that 
$$ |M(u)-M(v)| \lesssim {|u-v|}^{\alpha} \int {|\widehat{\mu}(\xi)|}^2 {|\xi|}^{-\frac{d-1}{2}+\alpha} d\xi 
\lesssim {|u-v|}^{-\alpha},$$ where the last step follows by Lemma \ref{energy}, and so the item iv) follows.

\end{document}